\documentclass[11pt]{article}

\usepackage{amsmath, amssymb, amsthm}
\usepackage{stmaryrd}
\usepackage[top = 1in, bottom = 1in, left = 1in, right = 1in]{geometry}
\newtheorem{thm}{Theorem}[section]
\newtheorem{lem}[thm]{Lemma}

\theoremstyle{remark}
\newtheorem{remark}[thm]{Remark}

\newcommand{\EV}{\operatorname*{\mathbb{E}}}
\newcommand{\RN}{\mathbb{R}}

\newcommand{\paw}{\operatorname{paw}}
\newcommand{\ind}{\operatorname{ind}}
\newcommand{\inj}{\operatorname{inj}}
\newcommand{\Ind}{\text{Ind}}

\newcommand{\Tr}{\text{Tr}}
\newcommand{\Span}{\text{Span}}
\newcommand{\A}{\mathcal{A}}
\newcommand{\Fl}{\mathcal{F}}
\newcommand{\K}{\mathcal{K}}
\newcommand{\W}{\mathcal{W}}
\newcommand{\X}{\mathbf{X}}
\newcommand{\F}{\mathbf{F}}
\newcommand{\Y}{\mathbf{Y}}
\newcommand{\M}{\mathbf{M}}

\title{The Inducibility of Graphs on Four Vertices}
\author{James Hirst \\ School of Computer Science, McGill University, Montreal, Canada \\ james.hirst@mail.mcgill.ca}
\date{}

\begin{document}
\maketitle

\begin{abstract}
We consider the problem of determining the maximum induced density of a graph $H$ in any graph on $n$ vertices. The limit of this density as $n$ tends to infinity is called the inducibility of $H$. The exact value of this quantity is known only for a handful of small graphs and a specific set of complete multipartite graphs. Answering questions of Brown-Sidorenko and Exoo we determine the inducibility of $K_{1,1,2}$ and the paw graph.  The proof is obtained using semi-definite programming techniques based on a modern language of extremal graph theory, which we describe in full detail in an accessible setting.

\end{abstract}

\section{Introduction}
Following the notation of \cite{Exo86}, for graphs $H,G$ with $|V(H)| = k$, $|V(G)| = n$ we define $I(H;G)$ as the number of induced subgraphs of $G$ isomorphic to $H$. Since we are not concerned with the size of $G$, we normalize $I(H;G)$ in a suitable way by defining $i(H;G) = I(H;G)/\binom{n}{k}$. We are now able to define the \textit{inducibility} of a graph $H$ as
\[ i(H) = \lim_{n \rightarrow \infty} \max_{\substack{ G \\ |V(G)| = n}} \ i(H;G). \]
%
It is not difficult to see that this limit always exists. The very recent work \cite{HHN} gives some strong asymptotic results for a large class of graphs, but the problem of determining $i(H)$ appears to be non-trivial even in some cases when $H$ is a very small graph. Let $\overline{H}$ denote the complement graph of $H$. Note that we have that $i(H) = i(\overline{H})$, so that we need only to consider one graph in a complementary pair. The inducibility of all graphs on less than four vertices is known: The only such non-trivial graphs are $P_3$ and its complement, and it is shown in \cite{Exo86} that $i(P_3) = 3/4$. There are 11 non-isomorphic graphs on four vertices: five complementary pairs and one self-complementary graph. The case of $K_{4}$, the complete graph on four vertices is trivial, as for any complete graph $K_{n}$ we have $i(K_{n}) = i(\overline{K_{n}}) = 1$.

In \cite{Exo86}, Exoo gave upper and lower bounds on $i(H)$ for the remaining 4-vertex graphs. Since then the inducibility of some of these graphs has been determined. In \cite{BNT86}, Bollob{\'a}s \textit{et al}. studied the inducibility of balanced complete bipartite graphs $K_{r,r}$, and proved that the graph which maximizes $i(K_{r,r};G)$ is again a balanced complete bipartite graph $K_{n,n}$. This resolved the case of $K_{2,2}$ and its complement. Later in \cite{BS94}, Brown and Sidorenko proved that if $H$ is a complete bipartite graph, then the graph which maximizes $i(H;G)$ can also be chosen to be complete bipartite, resolving $K_{1,3}$ and its complement. They also gave a generalization to complete multi-partite graphs, along with conditions under which the exact graph $G$ is known. However, the inducibility of the complete tri-partite $K_{1,1,2}$ is not classified by these conditions.

There remain three 4-vertex graphs (considering one from each complementary pair) for which the exact value of $i(H)$ is unresolved: $K_{1,1,2}$, the paw graph (the graph obtained from a triangle by adding a pendant edge) $H_{\paw}$, and the path on four vertices $P_{4}$. Brown and Sidorenko mentioned in \cite{BS94} that the best construction they know for $K_{1,1,2}$ is a balanced complete multipartite graph with 5 parts (a ``5-equipartite graph'' by their convention). Making use of the recent theory of flag algebras from \cite{Raz07} and semi-definite programming techniques used in \cite{Raz10,HKN,Grzesik,HHKNR1,HHKNR2} among other papers, we determine $i(K_{1,1,2})$, and $i(H_{\paw})$, and show that the construction of Brown and Sidorenko is in fact best possible.

A detailed account of the semi-definite method in extremal graph theory is not easily accessible. The description given in \cite{Raz07} is difficult, as it is presented in a very general language, and other accounts are often given in terms of hypergraphs (\cite{Raz10}), so we will give a full development of the methods by which we obtained our results. Our notation differs slightly from that of the flag algebras introduced by Razborov in \cite{Raz07}, so we present in the next section the basic definitions which will be needed. Following this, we present the semi-definite program which does the bulk work in obtaining our results, and show how to interpret the output of the solver. Finally, we give techniques, some known and some new, for obtaining the best possible proof.

\section{Preliminaries}
\subsection{Homomorphism Densities}
For two graphs $H$ and $G$, a graph \textit{homomorphism} from $H$ to $G$ is a map $\varphi:V(H) \rightarrow V(G)$ which preserves adjacency, i.e. $\varphi$ is such that $(u,v) \in E(H) \implies (\varphi(u),\varphi(v)) \in E(G)$ for every $u,v \in V(H)$. Now we can define $t(H;G)$ as the probability that a uniformly chosen map $\varphi: V(H) \rightarrow V(G)$ is a homomorphism. The quantity $t(H;G)$ is called the homomorphism density of $H$ in $G$. We similarly define $t_{\ind}(H;G)$ with the additional condition that the homomorphisms should also preserve non-adjacency.

While $t(H;G)$ is interesting in its own right, extremal graph theory more often studies the quantity $t^{\inj}(H;G)$, which is defined as the probability that a uniformly chosen \emph{injective} map $\varphi$ is a homomorphism. However, the following lemma due to Lov{\'a}sz and Szegedy in \cite{LS06} shows that the two are close up to an error term of $o(1)$.
\begin{lem}
\label{lem:tAsym}
For every two graphs $H$, and $G$,
\[ |t(H;G) - t^{\inj}(H;G)| \leq \frac{1}{|V(G)|}\binom{|V(H)|}{2} = o_{|V(G)| \rightarrow \infty}(1). \]
\end{lem}
For our purposes, Lemma \ref{lem:tAsym} will allow us to restrict our attention to the more nicely behaved function $t(H;\cdot)$. Now we give an alternative, but useful definition of $t(H;G)$ and $t_{\ind}(H;G)$ in terms of $A_{G}$, the adjacency matrix of $G$. For graphs $H$ and $G$, let $\{x_{u} \ | \ u \in V(H)\}$ be independent random variables which take values from $V(G)$ uniformly. Then we have
\begin{equation}
\label{eq:tA}
t(H;G) = \EV \left[\prod_{\substack{(u,v)\in E(H)}} A_{G}(x_{u},x_{v})\right],
\end{equation}
and
\begin{equation}
\label{eq:tAind}
t_{\ind}(H;G) = \EV \left[\prod_{(u,v) \in E(H)} A_{G}(x_u,x_v) \prod_{(u,v) \notin E(H)} (1 - A_{G}(x_{u},x_{v}))\right].
\end{equation}
The two functions are related by
\begin{equation}
\label{eq:Ind}
t(H;\cdot) = \sum_{\substack{F\supseteq H \\ V(F) = V(H)}} t_{\ind}(F;\cdot),
\end{equation}
and a M\"{o}bius inversion formula
\begin{equation}
\label{eq:Inv}
t_{\ind}(H;\cdot) = \sum_{\substack{F \supseteq H \\ V(F) = V(H)}} (-1)^{|E(F) \backslash E(H)|}t(F;\cdot).
\end{equation}

\subsection{Graphons}
A sequence of graphs $\{G_{i}\}_{i=1}^{\infty}$ is called \textit{convergent} if, for every graph $H$, the sequence $\{t(H;G_{i})\}_{i=1}^{\infty}$ converges.
It is not difficult to construct convergent sequences $\{G_i\}_{i=1}^\infty$ such that their limits cannot be recognized as graphs, i.e. there is
no graph $G$ with $\lim_{i \rightarrow \infty} t(H;G_i) = t(H;G)$ for every $H$. For this reason, the extremal solution to a problem is often stated as a sequence of graphs rather than a single graph. However, there is a ``soft analytic'' approach which enables us to avoid working with sequences. We can extend the space of graphs, and represent the limits of convergent sequences of graphs as an object in this extended space. It is shown in~\cite{LS06} that the limit of a convergent graph sequence can be represented as a measurable, symmetric function $w:[0,1]^{2} \rightarrow [0,1]$. Let $\W$ denote the set of all such functions. The elements of $\W$ are called \textit{graphons}. We can extend the definition of the functions $t(H;\cdot)$ and $t_{\ind}(H;\cdot)$ to graphons. Let $\{x_{u} \ | \ u \in V(H)\}$ be independent random variables which take values uniformly from $[0,1]$. Then for $w \in \W$, we define
\begin{equation}
\label{eq:tG}
t(H;w) = \EV \left[\prod_{\substack{(u,v)\in E(H)}} w(x_{u},x_{v})\right],
\end{equation}
and
\begin{equation}
\label{eq:tGind}
t_{\ind}(H;w) = \EV \left[\prod_{(u,v) \in E(H)} w(x_{u},x_{v}) \prod_{(u,v) \notin E(H)} (1 - w(x_{u},x_{v}))\right].
\end{equation}
Let $G$ be a graph with $|V(G)| = n$, and define a graphon $w_{G}$ as follows. Let $w_{G}(x,y) = A_{G}(\lceil xn \rceil, \lceil yn \rceil)$ for $(x,y) \in (0,1]^{2}$, and $w_{G} = 0$ otherwise. Comparing (\ref{eq:tA}) and (\ref{eq:tG}), it is easy to see that we have $t(H;G) = t(H;w_{G})$ (and $t_{\ind}(H;G) = t_{\ind}(H;w_{G})$ similarly) for every $H,G$. In this sense graphons are a natural extension of finite graphs. Furthering this notion, Lov{\'a}sz and Szegedy show in \cite{LS06} that the space of graphons $\W$ is complete in the sense that every convergent sequence of graphons converges to a graphon, and also that the set $\{w_{G} \ | \ G \text{ is a finite graph} \}$ is dense in $\W$.
That is, if $\{w_i\}_{i \in \mathbb{N}}$ is a convergent sequence of graphons, then there exists a graphon $w$ such that $\lim_{i \to \infty} t(H;w_i)=t(H;w)$ for every finite graph $H$, and furthermore for every graphon $w$, there exists a sequence of finite graphs $\{G_i\}_{i \in \mathbb{N}}$ such that $\lim_{i \to \infty} t(H;w_{G_i})=t(H;w)$ for every graph $H$.

Graphons allow us to give a simpler analogue of $i(H)$ in the language of homomorphisms. For a graph $H$, let $\Gamma(H)$ denote its group of automorphisms. Then we define
\begin{equation}
\label{eq:tind}
t_{\ind}(H) = \max_{w \in \W} \ t_{\ind}(H;w) = \frac{|\Gamma(H)|}{|V(H)|!} i(H).
\end{equation}
%
Note that the fact that the maximum in (\ref{eq:tind}) is attained follows from the completeness of the space of graphons and the compactness of $[0,1]^\mathbb{N}$.

\subsection{Quantum Graphs}
A $k$\textit{-partially labeled graph} is a graph in which $k$ of the vertices have been labeled with distinct natural numbers $1,\ldots,k$. We extend the definition of $t(\cdot;w)$ for partially labeled graphs. For a $k$-partially labeled graph $H$, and a map $\phi:[k] \rightarrow [0,1]$, we define $t(H,\phi;w)$ as $t(H;w)$ in (\ref{eq:tG}), conditioned on the event that $x_{i} = \phi(i)$, $1 \leq i \leq k$.

The function $\llbracket  \cdot \rrbracket$ maps $k$-partially labeled graphs to graphs by unlabeling the labeled vertices. This function can be seen as an averaging map, as for a $k$-partially labeled graph $H$, we have
\begin{equation}
\label{eq:unlab}
t(\llbracket  H \rrbracket;w) = \EV_{\phi} t(H,\phi;w),
\end{equation}
where $\phi:[k] \rightarrow [0,1]$ is a uniformly chosen map. We define the product of $k$-partially labeled graphs $H_{1}$ and $H_{2}$, denoted $H_{1} \cdot H_{2}$, by taking the disjoint union of $H_{1}$ and $H_{2}$ as graphs and then identifying the labeled vertices, reducing multiple edges. Note that by this definition we have the property that
\begin{equation}
\label{eq:algebraMorphism}
t(F \cdot G, \phi;w) = t(F,\phi;w)t(G,\phi;w),
\end{equation}
for $k$-partially labeled graphs $F,G$ and any compatible map $\phi$.

Let $\Fl_{k}$ denote the set of all $k$-partially labeled graphs. Together with the above product, $\Fl_{k}$ has the structure of a commutative semigroup. Now let $\RN [\Fl_{k}]$ denote the semigroup algebra of $\Fl_{k}$ over $\RN$. That is, elements of $\RN [\Fl_{k}]$ are formal finite linear combinations of elements of $\Fl_{k}$, and the product on $\RN [\Fl_{k}]$ is obtained from the product on $\Fl_{k}$ by distributivity over addition. We extend the definition of the functions $t(\cdot,\phi;w)$ and $\llbracket \cdot \rrbracket$ to elements of $\RN [\Fl_{k}]$ in the natural way (by linearity).
It follows from (\ref{eq:algebraMorphism}) that for every fixed graphon $w$ and map $\phi:[k] \rightarrow [0,1]$, the function $t(\cdot,\phi;w):\RN [\Fl_{k}] \to \mathbb{R}$ is an algebra homomorphism. Let $\K_k$ denote the intersection of the kernels of all these algebra homomorphisms. In other words,
$$\K_{k} = \{ f \in \RN [\Fl_{k}] \ | \ t(f, \phi;w) = 0 \quad \forall \phi,w \}.$$
Then $\K_k$ is a subalgebra and we obtain the quotient algebra $\A^{k}$ by $\A^{k} = \RN [\Fl_{k}] / \K_{k}$. Elements of $\A^k$ are called $k$-\textit{quantum graphs} (we will refer to 0-quantum graphs as simply quantum graphs).
Note that for every $k$-partially labeled graph $H \in \A^k$, we have
\begin{equation}
\label{eq:lift}
H = H \cup K_1.
\end{equation}

For a $k$-partially labeled graph $H$, the $k$-quantum graph $\Ind(H)$ is defined as
\[ \Ind(H) := \sum_{\substack{F \supseteq H \\ V(F) = V(H)}} (-1)^{|E(F) \backslash E(H)|}F. \]
Note that it follows from (\ref{eq:Inv}) that $t_{\ind}(H,\phi;w) = t(\Ind(H),\phi;w)$ for every $\phi,w$.

A partially labeled graph in which all the vertices are labeled is called a \textit{type}. For a fixed type $\sigma$, let $H$ be a partially labeled graph for which the labeled vertices induce a subgraph isomorphic (by a label preserving isomorphism) to $\sigma$. Then $\Ind(H) \in \A^{k}$ is called a $\sigma$-\textit{flag}. Let $\Ind(F),\Ind(H)$ be $\sigma_{1},\sigma_{2}$-flags, respectively, with $\sigma_{1} \neq \sigma_{2}$. Then $\Ind(F) \cdot \Ind(H) = 0$, since for any $\phi,w$, at least one of $t(\Ind(F),\phi;w)$, $t(\Ind(H),\phi;w)$ is zero.

For every type $\sigma$ on $k \geq 1$ vertices, let $\A^{\sigma}$ be the subalgebra of $\A^{k}$ spanned by the set $\{ f \in \A^{k} \ | \ f \text{ is a } \sigma \text{-flag} \}$. If $\sigma_{1},\ldots,\sigma_{r}$ are all the non-isomorphic types on $k$ vertices, then we have the orthogonal decomposition $\A^{k} = \oplus_{i=1}^{r} \A^{\sigma_{i}}$. In other words, for every $f \in \A^{\sigma_{i}}$ and $g \in \A^{\sigma_{j}}$ with $i \neq j$, we have $f \cdot g = 0$.

\begin{remark}
The reader should be warned that our notation here slightly differs from that of Razborov in \cite{Raz07}. A $\sigma$-flag $H$ in Razborov's language corresponds to $(u!/|\Gamma(H)|)\Ind(H)$ in our language, where $u$ is the number of unlabeled vertices in $H$, and $\Gamma(H)$ refers here to the group of label preserving automorphisms of $H$. We believe that our notation has the advantage that the product of the algebra $A^\sigma$, and the operator $\llbracket \cdot \rrbracket$ have easier descriptions.
\end{remark}

The crux of our results relies on the following easy lemma.
\begin{lem}
\label{lem:pos}
Let $f$ be a $k$-quantum graph, then
\[ t(\llbracket  f^{2} \rrbracket;w) \geq 0 \]
for all graphons $w$.
\begin{proof}
Fix a map $\phi:[k] \rightarrow [0,1]$, and let $f = \sum_{i=1}^{n} \alpha_{i} H_{i}$. Then we have
\[ t(f^{2},\phi;w) = t\left( \left( \sum_{i=1}^{n} \alpha_{i} H_{i} \right)^{2}, \phi; w\right) = t\left( \sum_{i=1}^{n} \sum_{j=1}^{n} \alpha_{i} \alpha_{j} H_{i} \cdot H_{j}, \phi; w\right) \]
\[ = \sum_{i=1}^{n} \sum_{j=1}^{n} \alpha_{i} \alpha_{j} t(H_{i}, \phi; w) t(H_{j},\phi;w) = \left( \sum_{i=1}^{n} \alpha_{i} t(H_{i},\phi;w) \right)^{2} \geq 0. \]
Now our desired result follows from (\ref{eq:unlab}), which shows that the function $ \llbracket \cdot \rrbracket $ preserves positivity.
\end{proof}
\end{lem}
\begin{remark}
\label{rem:PSD}
It follows from Lemma $\ref{lem:pos}$ that for any $k$-quantum graphs $f_{1},\ldots,f_{n}$, and any graphon $w$, the $n \times n$ matrix $\M$ defined by $\M(i,j) = t(\llbracket f_{i}\cdot f_{j}\rrbracket;w)$ is positive semi-definite (PSD). To see this, note that we can write
\[ z^{T} M z = \sum_{i,j} z_{i}z_{j}t(\llbracket f_{i} \cdot f_{j}\rrbracket;w) = t\left( \left\llbracket  \left(\sum_{i} z_{i} f_{i}\right)^{2}\right\rrbracket; w\right) \geq 0. \]
This is the basis for the so called semi-definite method in extremal graph theory.
\end{remark}

\section{The Semi-Definite Method}
\subsection{The Semi-Definite Program and its Dual \label{sec:SDP}}
It follows from reflection positivity characterizations \cite{FLS07, LS06, LS09} that every asymptotic algebraic inequality between subgraph densities follows from the positive semi-definiteness of a certain infinite matrix. If this matrix were finite, then it would be possible to solve every such inequality efficiently using semi-definite programming. However, since the matrix is infinite, in practice one can only consider a finite sub-matrix and hope that the inequality still holds under the weaker condition that the sub-matrix is positive semi-definite. We call this approach the semi-definite method in extremal graph theory. First let us recall the standard form of a semi-definite program (SDP) and its dual.

The primal SDP formulation we will make use of is the following
\begin{equation}
\label{eq:primal}
\begin{array}{ll}
\text{minimize } & c \cdot x \\
\text{subject to } & \X = \sum_{i=1}^{m} \F_{i} x_{i} - \F_{0} \succeq 0
\end{array}
\end{equation}
and the corresponding dual program is
\begin{equation}
\label{eq:dual}
\begin{array}{ll}
\text{maximize } & \F_{0} \bullet \Y \\
\text{subject to } & \F_{i} \bullet \Y = c_{i} \quad 1 \leq i \leq m, \quad \Y \succeq 0
\end{array}
\end{equation}
where $c \in \RN^{m}$, $\F_{i} \in \RN^{n \times n}$ are given, and $x \in \RN^{m}, \Y \in \RN^{n \times n}$ are variables. Note that $c \cdot x$ refers to the standard dot product, while $\F_{0} \bullet \Y$ refers to the inner product on $\RN^{n\times n}$ defined by $A \bullet B = \Tr(A^{T}B) = \sum_{i,j}A_{i,j}B_{i,j}$.
The duality here refers to the fact (see e.g.~\cite{MR1952986}) that feasible solutions to (\ref{eq:primal}) give upper bounds for feasible solutions to (\ref{eq:dual}) and vice-versa.

Now fix a quantum graph $f = \sum_{i=1}^{m} \alpha_{i} H_{i}$, and suppose we want to find
\begin{equation}
\label{eq:obj}
\min_{w \in \W} t(f;w).
\end{equation}
%
Many statements in asymptotic extremal graph theory can be converted to this form. For example, the celebrated Goodman bound~\cite{MR0107610} says that
$\min_{w \in \W} t(K_3-2(K_2 \cup K_2)+K_2;w)  \ge 0$. 

Choose two parameters $0 \leq L \leq N$, where $N$ is large enough so that  
$$f \in \Span \{H :  \mbox{$H$ is an $N$-vertex graph}\} \subseteq \mathcal{A}^0.$$
For every $1 \leq k \leq L$, let $\{f_{1}^{k},\ldots,f_{m_{k}}^{k}\}$ be a set of independent $k$-quantum graphs on $\lfloor (N+k)/2 \rfloor$ vertices, so that by (\ref{eq:lift}) we have
$$f_{i}^{k} \cdot f_{j}^{k} \in \Span \{H :  \mbox{$H$ is an $N$-vertex $k$-partially labeled graph}\} \subseteq \mathcal{A}^k.$$

By (\ref{eq:lift}), in the expansion $f = \sum_{i=1}^m \alpha_i H_i$, we can assume that $H_1,\ldots,H_m$ are exactly all the non-isomorphic graphs on $N$ vertices. Consider the following semi-definite program:
\begin{equation}
\label{eq:opt}
\begin{array}{ll}
\text{minimize } & \sum_{i=1}^{m} \alpha_{i} x_{H_{i}} \\
\text{subject to } & x_{K_{1}} = 1; \\
 & \forall k \leq L \text{, the matrix $\M$ defined by }  \M(u,v) = \sum_{i=1}^m \beta_i x_{H_i} \text{ is PSD,} \\
 & \text{where }\sum_{i=1}^m \beta_i H_i := \llbracket f_{u}^{k} \cdot f_{v}^{k}\rrbracket.
\end{array}
\end{equation}

One should interpret the values $x_{H}$ in $(\ref{eq:opt})$ as $t(H;w)$, for some $w$. Note then that Remark $\ref{rem:PSD}$ shows that these semi-definiteness constraints are in fact necessary. By this interpretation, the solution to $(\ref{eq:opt})$ will be a lower bound for $(\ref{eq:obj})$. By choosing $L$ and $N$ large enough, this lower bound will often be quite good, and in some cases best possible. This is not the case for any choice of $f$, however. In \cite{HN11}, it is shown that there exists quantum graphs $f$ for which no choice of $L$ and $N$ will give a sharp bound. In fact, they show that determining the value of $(\ref{eq:obj})$ for general quantum graphs $f$ is undecidable.

It is not too hard to see how to convert $(\ref{eq:opt})$ to the form $(\ref{eq:primal})$. To implement the constraints as a semi-definiteness condition, we make use of the fact that a block diagonal matrix is PSD if and only if all of its diagonal blocks are. Each constraint is implemented using one block. 
For the first $L$ blocks, the $uv$-th entry in the $k$-th block of $\F_{i_0}$ is given by $\beta_{i_0}$, where $\llbracket f_{u}^{k} \cdot f_{v}^{k}\rrbracket =  \sum_{i=1}^m \beta_i H_i$. We use the last  block (i.e. the ($L+1$)-th block) to implement the first constraint: We want $x_{K_{1}} - 1 \geq 0$ and $1 - x_{K_{1}} \geq 0$, which are implemented making use of the constant matrix $\F_{0}$.

\subsection{Extracting a Lower Bound}
In Section~\ref{sec:SDP}, we saw  how (\ref{eq:opt}) can be converted to the form (\ref{eq:primal}). Now suppose the optimal solution to the dual problem (\ref{eq:dual}) is $\alpha$. Then our above discussion shows that $\alpha$ is a lower bound for $t(f;w)$. The matrix $\Y$ from the dual problem encodes a proof of this, which we can extract. The constraint on $\Y \succeq 0$ is that $\F_{i} \bullet \Y = \alpha_{i}$. So we have
\[ \sum_{i=1}^{m} ( \F_{i} \bullet \Y ) H_{i} = \sum_{i=1}^{m} \alpha_{i} H_{i} = f. \]
%
For every $1 \le k \le L$, let $\Y^{(k)}$  denote the $k$-th block of $\Y$. That is the block that corresponds to the positive semi-definiteness constraint created by $\{f_{1}^{k},\ldots,f_{m_{k}}^{k}\}$ in (\ref{eq:opt}). Similarly let $\F_i^{(k)}$ denote the $k$-th block of $\F_i$ for every $1 \le i \le m$ and every $1 \le k \le L$.

By splitting the inner product between the last block (the $x_{K_{1}} = 1$ constraint) and the rest of the blocks, it follows from $\F_{0} \bullet Y = \alpha$ that
\[ \sum_{i=1}^{m} \sum_{k=1}^L \left( \sum_{u=1}^{m_k} \sum_{v=1}^{m_k} \F_{i}^{(k)}(u,v) \Y^{(k)}(u,v)\right) H_{i} = f - \alpha K_{1}=f-\alpha. \]
Here we use $\alpha$ to denote the quantum graph $\alpha K_{1}$ as a slight abuse of notation, based on the fact that $t(\alpha K_{1};w) = \alpha$ for all graphons $w$.
Exchanging the order of summation, and using the definition of $\F_\ell$, this reduces to
\[f - \alpha = \sum_{k=1}^L \sum_{u=1}^{m_k} \sum_{v=1}^{m_k} \Y^{(k)}(u,v)\llbracket f^k_{u} \cdot f^k_{v}\rrbracket = \sum_{k=1}^L \llbracket z_k^{T}\Y^{(k)} z_k\rrbracket, \]
where $z_k$ is the vector $(f^k_{1},\ldots,f^k_{m_k})^T$. 
\begin{remark}
\label{rem:PSDform}
For any matrix $A \succeq 0$ and vector of $k$-quantum graphs $z$, the quantum graph $\llbracket z^{T} A z \rrbracket$ is trivially positive. To see this, notice that for any map $\phi:[k] \rightarrow [0,1]$ and a graphon $w$, we can write
\[ t(z^{T} A z,\phi;w) = \sum_{i,j} A(i,j) t(z_{i},\phi;w) t(z_{j},\phi;w) = v^{T} A v \geq 0, \]
where $v$ is defined by $v_{i} = t(z_{i},\phi;w)$. Now the positivity of $\llbracket z^{T} A z \rrbracket$ follows as in Lemma \ref{lem:pos}.
\end{remark}
In particular, Remark \ref{rem:PSDform} shows that the quantum graph $f - \alpha$ is positive, i.e. $t(f;w) \geq \alpha$ for every graphon $w$.
\begin{remark}
\label{rem:eigenform}
Since $\Y$ is positive semi-definite, every block $\Y^{(k)}$ is also positive semi-definite. Performing an eigenvalue decomposition on each $\Y^{(k)}$, and noting that $\Y^{(k)} \succeq 0$ has non-negative eigenvalues, we can get the useful form
\[ f - \alpha = \sum_{i=1}^{r} \lambda_{i} \llbracket g_{i}^2 \rrbracket, \]
where each $g_{i} \in \Span(\{f_{1}^{k},\ldots,f_{m_{k}}^{k}\})$, for some $k$.
\end{remark}

\subsection{Some Remarks}
From Remark \ref{rem:eigenform}, we can see that the choice of sets $\{f_{1}^{k},\ldots,f_{m_{k}}^{k}\}$ is key in obtaining a proof of our desired bound. 
%
%
To have the best chance of finding a proof, we should choose $f_{1}^{k},\ldots,f_{m_{k}}^{k}$ so that they span the entire algebra $\A^{k}$. One simple choice is to use the set of all $k$-partially labeled graphs on $\lfloor (N+k)/2 \rfloor$ vertices. However, the resulting matrix $\Y$ will be dense, and it is in our interests to obtain as small and simple a proof as possible. As a first measure, note that the size of $\Y$ increases exponentially with the parameters $L$ and $N$, so we should choose these to be minimal while still getting the desired bound.

The set $$\left\{ \Ind(F) : \mbox{$F$ is a $\lfloor (N+k)/2 \rfloor$-vetex $k$-partially labeled graph}\right\}$$ makes a particularly good basis for $\A^{k}$. It follows easily from (\ref{eq:Ind}) and (\ref{eq:Inv}) that this is indeed a basis for $\A^{k}$. We can partition this set according to the type obtained by restricting each $F$ to its labeled vertices. Then the partition corresponding to each type $\sigma$ is a basis for $\A^{\sigma}$. Because of the orthogonality properties these subspaces enjoy, by this choice of basis each diagonal block of the matrix $\Y$ reduces to a number of smaller blocks: one for each type on $k$ vertices.

We should also note here that not even all these blocks are required. By applying the eigenvalue decomposition used in Remark \ref{rem:eigenform} to the new, sparser $\Y$ matrix, we will get that $f - \alpha = \sum_{i=1}^{r} \lambda_{i} \llbracket g_{i}^2 \rrbracket$, where now each $g_i \in \A^{\sigma}$ for some type $\sigma$. If $\sigma_{1} \neq \sigma_{2}$ are two types such that $\llbracket \sigma_{1} \rrbracket = \llbracket \sigma_{2} \rrbracket$, then for every $g \in \A^{\sigma_{1}}$, there exists some $f \in \A^{\sigma_{2}}$ such that $\llbracket g \rrbracket = \llbracket f \rrbracket$. Thus we may restrict our search space to include only one of $\A^{\sigma_1},\A^{\sigma_2}$. In particular, we only need a block in $\Y$ for each type taken up to (not necessarily label preserving) isomorphism.

The structure of $A^{\sigma}$ is studied in detail in \cite{Raz07}, and it provides a decomposition of $A^{\sigma}$ as a direct sum of two smaller algebras. The group $\Gamma(\llbracket \sigma \rrbracket)$ acts on $\A^{\sigma}$ by permuting the labels. Then it can be shown that we can decompose $\A^{\sigma}$ orthogonally as $\A^{\sigma,+} \oplus \A^{\sigma,-}$ (the ``invariant'' and ``anti-invariant'' parts) where
\[ \A^{\sigma,+} = \{ f \in \A^{\sigma} \ | \ \pi(f) = f \quad \forall \pi \in \Gamma(\llbracket \sigma \rrbracket)\}, \]
and
\[ \A^{\sigma,-} = \{ f \in \A^{\sigma} \ | \ \sum_{\pi \in \Gamma(\llbracket \sigma \rrbracket)} \pi(f) = 0 \}, \]
so that if $f \in \A^{\sigma,+}$, and $g \in \A^{\sigma,-}$, then $f \cdot g = 0$. So for each $\sigma$, choose a basis for $A^{\sigma}$ according to this decomposition, i.e. as the union of a basis for $\A^{\sigma,+}$ and a basis for $\A^{\sigma,-}$. Then each block of $\Y$ further reduces to two smaller blocks.

If we have a set $\W_{0}$ of ``conjectured extremal graphons'', i.e. graphons such that $t(f;w_{0}) = \alpha$, then we can also make use of the following. Let $w_{0} \in \W_{0}$. Then if
\[ \sum_{i=1}^{r} \llbracket g_{i}^2\rrbracket = f - \alpha, \]
it follows that
\[ \sum_{i=1}^{r} t(\llbracket g_{i}^{2}\rrbracket;w_{0}) = 0. \]
Thus for all $1 \leq i \leq r$,  we must have $$\EV_\phi t( g_{i}, \phi ;w_{0})^2= \EV_\phi t( g_{i}^2, \phi ;w_{0}) =  t(\llbracket g_{i}^2\rrbracket ;w_{0}) = 0,$$
which shows that  $\EV_\phi \left|t( g_{i}, \phi ;w_{0})\right|=0$.  Define $\Delta^{\sigma} = \{f \in A^{\sigma} \ | \ \EV_\phi \left|t(f, \phi ;w_{0})\right| = 0 \quad \forall w_{0} \in \W_{0} \}$. Then we can reduce the size of each block of $\Y$ by restricting ourselves to a basis of $\A^{\sigma} \cap \Delta^{\sigma} = (\Delta^{\sigma} \cap \A^{\sigma,+}) \oplus (\Delta^{\sigma} \cap \A^{\sigma,-})$.

A final remark is in regards to the floating point nature of the SDP solver. The matrix $\Y$ output by the SDP corresponds to a proof of the desired bound only to a specified degree of floating point error. While some of the entries may have obvious rational closed forms (computer algebra systems are quite good at finding these), it is likely that the solver will have more degrees of freedom than required. Thus the entries will often be somewhat arbitrary floating point numbers. To resolve this, we can first attempt to restrict our bases even further. We are not aware of any obvious way to do this, and the best results will come from repeated experiments. Once we have a minimal size proof, to resolve any remaining slackness in the solver, we can introduce further constraints to fix certain values of $\Y$ to close rational approximations. By fixing only a few values, the remaining entries will often be uniquely determined rational numbers, and the solver does most of the work here.

\section{Results}
\subsection{Representing Graphs}
Before we give our results, we will need a notation that is able to concisely express a large number of partially labeled graphs. We will use an adjacency list, slightly modified for partially labeled graphs. Labeled vertices are identified by their labels, and other vertices are identified with a character. As an example, $\{1a, 1b, 2a, 2b, ab\}_{4, 2}$ refers to a 4-vertex, 2-partially labeled graph (indicated in the trailing subscript) with five edges, and each pair corresponds to an edge.

\subsection{A Proof for the Paw Graph}
\begin{thm}
\label{thm:paw}
Let $H_{\paw}$ denote the four vertex graph obtained from a triangle by adding a single edge. Then we have
\[ t_{\ind}(H_{\paw}) = \frac{1}{32}. \]
\begin{proof}
We first give a graphon $w_{0}$ which achieves our desired lower bound. Let $w_{0} = 1 - w_{K_{2} \cup K_{2}}$. A simple counting argument now shows $t_{\ind}(H_{\paw};w_{0}) = 1/32$, giving the lower bound $t_{\ind}(H_{\paw}) \geq 1/32$.

To prove the upper bound, we must show that $t(\Ind(H_{\paw});w) \leq 1/32$ for all graphons $w$. Then we need to give positive semi-definite matrices $\Y_i$, and vectors of $k_{i}$-quantum graphs $z_i$ such that $\sum \llbracket  z_i^{T} \Y_i z_i \rrbracket = 1/32 - \Ind(H_{\paw})$. In fact if $H_1,\ldots,H_m$ are all graphs on $5$ vertices and $1/32 - \Ind(H_{\paw}) = \sum_{i=1}^m \alpha_i H_i$, then it suffices to have $\sum \llbracket  z_i^{T} \Y_i z_i \rrbracket = \sum_{i=1}^m \beta_i H_i$, where $\beta_i \le \alpha_i$ for every $1 \le i \le m$.

Define the positive semi-definite matrices $\Y_{1},\ldots,\Y_{6}$ as 
\[
\Y_{1} = \frac{1}{96}\begin{pmatrix} 4 & -7 & -2 & -5 & 4 \\
													 						 -7 & 59 & -38 & 33 & -7 \\
													 						 -2 & -38 & 44 & -18 & -2 \\
													 						 -5 & 33 & -18 & 19 & -5 \\
													 						 4 & -7 & -2 & -5 & 4 \end{pmatrix},
\quad
\Y_{2} = \frac{1}{1920}\begin{pmatrix} 80 & -275 & -70 \\
																 		  -275 & 1632 & -446 \\
																 		  -70 & -446 & 748 \end{pmatrix},
\]
\[
\Y_{3} = \frac{1}{192}\begin{pmatrix} 32 & -43 \\
																     -43 & 58 \end{pmatrix},
\quad
\Y_{4} = \frac{1}{960}\begin{pmatrix} 65 & -214 \\
																     -214 & 839 \end{pmatrix},
\]
\[
\Y_{5} = \frac{1}{12}\begin{pmatrix} 1 & 2 \\
																    2 & 4 \end{pmatrix},
\quad
\Y_{6} = \frac{1}{120}\begin{pmatrix} 24 & -13 \\
																     -13 & 10 \end{pmatrix},
\]
and the vectors $z_1,\ldots,z_6$  by
\[ z_1^{T} =  \Ind\left(\{ \}_{3, 1},\ \{1b\}_{3, 1},\ \{1b, ab\}_{3, 1},\ \{1a, 1b\}_{3, 1},\ \{1a, 1b, ab\}_{3, 1}\right), \]
\[ z_2^{T} = \Ind(\{3a\}_{4, 3} + \{2a\}_{4, 3},\ \{2a, 3a\}_{4, 3},\ \{1a, 2a\}_{4, 3}), \]
\[ z_3^{T} = \Ind(\{23\}_{4, 3},\ \{23, 2a, 3a\}_{4, 3} - \{1a, 23\}_{4, 3} - \{1a, 23, 3a\}_{4, 3}), \]
\[ z_4^{T} = \Ind(\{23, 3a\}_{4, 3} - \{23, 2a\}_{4, 3},\ \{1a, 23, 3a\}_{4, 3} - \{1a, 23, 2a\}_{4, 3}), \]
\[ z_5^{T} = \Ind(\{13, 23\}_{4, 3} - \{13, 23, 2a\}_{4, 3},\ \{13, 1a, 23, 2a\}_{4, 3} - \{13, 1a, 23, 2a, 3a\}_{4, 3}), \]
\[ z_6^{T} = \Ind(\{13, 23, 2a\}_{4, 3},\ \{13, 23, 2a, 3a\}_{4, 3} - \{13, 1a, 23, 3a\}_{4, 3}), \]
where as a slight abuse of notation, $\Ind(\cdot)$ applied to a row vector is defined by applying $\Ind(\cdot)$ to all its entries. It can be verified that $\Y_i$, and $z_i$ satisfy the conditions given above, so the desired upper bound holds, and the result follows.
\end{proof}
\end{thm}

\subsection{A Proof for $K_{1,1,2}$}
\begin{thm}
\label{thm:K112}
We have
\[ t_{\ind}(K_{1,1,2}) = \frac{12}{125} = 0.096. \]
\begin{proof}
The construction given by Brown and Sidorenko in \cite{BS94} corresponds to the graphon $w_{K_{5}}$, and a simple counting argument gives $t_{\ind}(K_{1,1,2};w_{K_{5}}) = 12/125$, giving us the lower bound.

The proof of the upper bound is as in the proof of Theorem \ref{thm:paw}, except that instead of $1/32 - \Ind(H_{\paw})$ we are dealing with $12/125 - \Ind(K_{1,1,2})$. We must also take $H_{1},\ldots,H_{m}$ (as defined in proof of Theorem \ref{thm:paw}) to be all the graphs on 7 vertices. The matrices  $\Y_{1},\ldots,\Y_{6}$ are defined as
\[
\Y_{1} = \frac{1}{180000}\begin{pmatrix} 17280 & -24912 & -79032 & 17916 & 27236 \\
													 						 -24912 & 70560 & 193926 & -71271 & -62329 \\
													 						 -79032 & 193926 & 587520 & -205728 & -185320 \\
													 						  17916 & -71271 & -205728 & 113760 & 52920 \\
													 						  27236 & -62329 & -185320 & 52920 & 62640 \end{pmatrix},
\]
\[
\Y_{2} = \frac{1}{7500}\begin{pmatrix} 270 & 2991 & -90 \\
																 		  2991 & 33150 & -997 \\
																 		  -90 & -997 & 30 \end{pmatrix},
\]
\[
\Y_{3} = \frac{1}{5000}\begin{pmatrix} 3075 & 5514 & 1483 & 3917 \\
													 						5514 & 14500 & 8915 & 12585 \\
													 						1483 & 8915 & 9320 & 9490 \\
													 						3917 & 12585 & 9490 & 11725 \end{pmatrix},
\quad
\Y_{4} = \frac{1}{100}\begin{pmatrix} 23 \end{pmatrix},
\]
\[
\Y_{5} = \frac{1}{300000}\begin{pmatrix} 95040 & -81360 & -68250 & -225285 & 43415 \\
																				-81360 & 417600 & 350340 & 749826 & -222822 \\
																				-68250 & 350340 & 294000 & 629130 & -186960 \\
																				-225285 & 749826 & 629130 & 1425600 & -400105 \\
																				43415 & -222822 & -186960 & -400105 & 107007 \end{pmatrix},
\]
\[
\Y_{6} = \frac{1}{1200000}\begin{pmatrix} 343560 & -310680 & -234812 & 53570 \\
																				-310680 & 190272 & 119819 & -34343 \\
																				-234812 & 119819 & 75456 & -21579 \\
																				 53570 & -34343 & -21579 & 9910 \end{pmatrix},
\]
and the vectors $z_1,\ldots,z_6$ are given by
\[ z_1^{T} =  \Ind(\{ \}_{4, 1},\ \{1c, ac, bc\}_{4, 1},\ \{1b, 1c, ab, ac\}_{4, 1},\ \{1a, 1b, 1c\}_{4, 1},\ \{1a, 1b, 1c, ac, bc\}_{4, 1}), \]
\[ z_2^{T} =  \Ind(\{1b, 2b, ab\}_{4, 2},\ \{1b, 2a, ab\}_{4, 2},\ \{1a, 1b, 2a, 2b, ab\}_{4, 2}), \]
\[ z_3^{T} =  \Ind(\{12\}_{4, 2},\ \{12, 2b\}_{4, 2} + \{12, 1b\}_{4, 2},\ \{12, 1b, 2b\}_{4, 2},\ \{12, 1b, 2a, 2b\}_{4, 2} + \{12, 1a, 1b, 2b\}_{4, 2}), \]
\[ z_4^{T} =  \Ind(\{12, 2b\}_{4, 2} - \{12, 1b\}_{4, 2}), \]
\[ z_5^{T} =  \Ind(\{13, 23, 3a, 3b\}_{5, 3},\ \{13, 1b, 23, 2b, 3a, 3b\}_{5, 3},\ \{13, 1b, 23, 2a, 2b, 3a\}_{5, 3} + \{13, 1a, 1b, 23, 2b, 3a\}_{5, 3},\]
\[ \indent \indent \indent \indent \indent \indent \indent \indent \{13, 1a, 1b, 23, 2a, 2b\}_{5, 3},\ \{13, 1a, 1b, 23, 2a, 2b, 3a, 3b\}_{5, 3}), \]
\[ z_6^{T} =  \Ind(\{12, 13, 1b, 23, 2a, 3a, 3b, ab\}_{5, 3} + \{12, 13, 1b, 23, 2a, 2b, 3a, ab\}_{5, 3} + \{12, 13, 1a, 1b, 23, 2b, 3a, ab\}_{5, 3}),\]
\[ \indent \indent \indent \indent \{12, 13, 1b, 23, 2a, 2b, 3a, 3b, ab\}_{5, 3} + \{12, 13, 1a, 1b, 23, 2b, 3a, 3b, ab\}_{5, 3}, \]
\[ \indent \indent \indent \indent \{12, 13, 1a, 1b, 23, 2a, 2b, 3b, ab\}_{5, 3},\ \{12, 13, 1a, 1b, 23, 2a, 2b, 3a, 3b, ab\}_{5, 3}). \]
\end{proof}
\end{thm}

\section{Conclusion and Open Problems}
The semi-definite method is a powerful technique for proving inequalities between subgraph densities, although it has its drawbacks. The proof we gave for the case of $K_{1,1,2}$ was obtained by efforts near the limit of our computing power. We had to run the semi-definite program with graphs on seven vertices (i.e. $N = 7$ in (\ref{eq:opt})). The solver took 2-3 hours to run on full precision, and the search space (and thus running time) grows exponentially with the number of vertices. Proofs by this method for larger graphs (7+ vertices) seem not to be feasible, but even some small graphs elude us. The only remaining unresolved graph on 4 vertices is the self-complementary $P_{4}$.

The case of $t_{\ind}(P_{4})$ is interesting. The upper bound given by Exoo in \cite{Exo86} was (converting to our language) $1/36 = 0.02\overline{77}$. Using the semi-definite method, we can improve this to a floating point bound of roughly $0.0172$, but we conjecture that this bound is not tight either. The best construction we know of for $P_{4}$ is still the one given by Exoo in $\cite{Exo86}$, involving a kind of blow-up of the Paley Graph $QR(17)$. The corresponding lower bound is $80/4877 \approx 0.0164$, so there is a clear gap between this construction and the best bounds we can obtain.

\section*{Acknowledgments}
My thanks are due to Hamed Hatami for introducing me to this problem, and providing many enlightening discussions regarding the details of this paper. Thanks are also due to NSERC, by whom my research was partially funded.

\bibliography{4graphs}

\begin{thebibliography}{10}

\bibitem{BNT86}
B{\'e}la Bollob{\'a}s, Chi{\^e} Nara, and Shun-ichi Tachibana.
\newblock The maximal number of induced complete bipartite graphs.
\newblock {\em Discrete Math.}, 62(3):271--275, 1986.

\bibitem{BS94}
Jason~I. Brown and Alexander Sidorenko.
\newblock The inducibility of complete bipartite graphs.
\newblock {\em J. Graph Theory}, 18(6):629--645, 1994.

\bibitem{Exo86}
Geoffrey Exoo.
\newblock Dense packings of induced subgraphs.
\newblock {\em Ars Combin.}, 22:5--10, 1986.

\bibitem{FLS07}
Michael Freedman, L{\'a}szl{\'o} Lov{\'a}sz, and Alexander Schrijver.
\newblock Reflection positivity, rank connectivity, and homomorphism of graphs.
\newblock {\em J. Amer. Math. Soc.}, 20(1):37--51 (electronic), 2007.

\bibitem{MR0107610}
A.~W. Goodman.
\newblock On sets of acquaintances and strangers at any party.
\newblock {\em Amer. Math. Monthly}, 66:778--783, 1959.

\bibitem{Grzesik}
Andrzej Grzesik.
\newblock On the maximum number of {$C_5$}'s in a triangle-free graph.
\newblock arXiv:1102.0962, 2011.

\bibitem{HHN}
Hamed Hatami, James Hirst, and Serguei Norine.
\newblock The inducibility of blow-up graphs.
\newblock arXiv:1108.5699, 2011.

\bibitem{HHKNR2}
Hamed Hatami, Jan Hladk{\'y}, Daniel Kr\'al', Serguei Norine, and Alexander
  Razborov.
\newblock Non-three-colorable common graphs exist.
\newblock arXiv:1105.0307, 2011.

\bibitem{HHKNR1}
Hamed Hatami, Jan Hladk{\'y}, Daniel Kr\'al', Serguei Norine, and Alexander
  Razborov.
\newblock On the number of pentagons in triangle-free graphs.
\newblock arXiv:1102.1634, 2011.

\bibitem{HN11}
Hamed Hatami and Serguei Norine.
\newblock Undecidability of linear inequalities in graph homomorphism
  densities.
\newblock {\em J. Amer. Math. Soc.}, 24(2):547--565, 2011.

\bibitem{HKN}
Jan Hladk\'y, Daniel Kr\'al', and Serguei Norine.
\newblock Counting flags in triangle-free digraphs.
\newblock arXiv:0908.2791, 2009.

\bibitem{MR1952986}
L{\'a}szl{\'o} Lov{\'a}sz.
\newblock Semidefinite programs and combinatorial optimization.
\newblock In {\em Recent advances in algorithms and combinatorics}, volume~11
  of {\em CMS Books Math./Ouvrages Math. SMC}, pages 137--194. Springer, New
  York, 2003.

\bibitem{LS06}
L{\'a}szl{\'o} Lov{\'a}sz and Bal{\'a}zs Szegedy.
\newblock Limits of dense graph sequences.
\newblock {\em J. Combin. Theory Ser. B}, 96(6):933--957, 2006.

\bibitem{LS09}
L{\'a}szl{\'o} Lov{\'a}sz and Bal{\'a}zs Szegedy.
\newblock Random graphons and a weak positivstellensatz for graphs.
\newblock {\em arXiv.org:0902.1327}, 2009.

\bibitem{Raz07}
Alexander~A. Razborov.
\newblock Flag algebras.
\newblock {\em J. Symbolic Logic}, 72(4):1239--1282, 2007.

\bibitem{Raz10}
Alexander~A. Razborov.
\newblock On 3-hypergraphs with forbidden 4-vertex configurations.
\newblock {\em SIAM J. Discrete Math.}, 24(3):946--963, 2010.

\end{thebibliography}

\end{document}